\documentclass{amsart}

\usepackage{url,amssymb,amsfonts, latexsym,amsmath,hhline,array,longtable}
\usepackage{pdfsync,color, comment,colortbl, mathrsfs,cite,graphicx}

\usepackage{xcolor}

\usepackage{ifpdf}
\ifpdf \usepackage[colorlinks=true, citecolor=blue,
linkcolor=blue, urlcolor=blue]{hyperref} \fi

\newcommand{\myx}{\boldsymbol{x}}

\newcommand{\cal}{\mathcal}

\def\avv#1{\textcolor{blue}{\bf #1}}

\newtheorem{formula}{}[section]
\newtheorem{definition}[formula]{Definition}
\newtheorem{corollary}[formula]{Corollary}
\newtheorem{remark}[formula]{Remark}
\newtheorem{lemma}[formula]{Lemma}
\newtheorem{theorem}[formula]{Theorem}
\newtheorem{example}[formula]{Example}

\def\thrm{\begin{theorem}}
\def\thrml#1{\begin{theorem}\label{#1}}
\def\ethrm{\end{theorem}}
\def\rmrk{\begin{remark}}
\def\rmrkl#1{\begin{remark}\label{#1}}
\def\ermrk{\end{remark}}
\def\xmpl{\begin{example}}
\def\xmpll#1{\begin{example}\label{#1}}
\def\exmpl{\end{example}}
\def\dfntn{\begin{definition}}
\def\dfntnl#1{\begin{definition}\label{#1}}
\def\edfntn{\end{definition}}
\def\nmrt{\begin{enumerate}}
\def\enmrt{\end{enumerate}}
\def\tm#1{\item[{\rm (#1)}]}
\def\qtnl#1{\begin{equation}\label{#1}}
\def\eqtn{\end{equation}}
\def\lmm{\begin{lemma}}
\def\lmml#1{\begin{lemma}\label{#1}}
\def\elmm{\end{lemma}}
\def\crllr{\begin{corollary}}
\def\crllrl#1{\begin{corollary}\label{#1}}
\def\ecrllr{\end{corollary}}
\def\css{\begin{cases}}
\def\ecss{\end{cases}}

\def\proof{\noindent{\bf Proof}.\ }

\def\cH{{\cal H}}

\def\cL{{\cal L}}

\DeclareMathOperator{\aut}{Aut}
\DeclareMathOperator{\alt}{Alt}
\DeclareMathOperator{\AGL}{AGL}
\DeclareMathOperator{\AGaL}{A{\rm \Gamma}L}

\DeclareMathOperator{\Fit}{Fit}

\DeclareMathOperator{\GaL}{{\rm \Gamma}L}
\DeclareMathOperator{\GF}{GF}
\DeclareMathOperator{\GL}{GL}

\DeclareMathOperator{\orb}{Orb}

\DeclareMathOperator{\Sp}{Sp}

\DeclareMathOperator{\sym}{Sym}

\DeclareMathOperator{\Z}{\mathbb{Z}}

\def\eprf{\hfill$\square$}

\def\qaq{\quad\text{and}\quad}

\def\m{^{(m)}}
\def\3{^{(3)}}
\def\ov{\overline}
\def\phmb#1{{\phantom{x}\hspace{-2mm}^{#1}}}

\def\wt{\widetilde}

\begin{document}

\title{The $3$-closure of a solvable permutation group is solvable}
\author{E.\,A.\ O'Brien}
\address{Department of Mathematics, The University of Auckland,
Auckland, New Zealand}
\email{e.obrien@auckland.ac.nz}
\author{I.\ Ponomarenko}
\address{Steklov Institute of Mathematics at St.\ Petersburg,  Russia;\newline
and Sobolev Institute of Mathematics, Novosibirsk, Russia}
\email{inp@pdmi.ras.ru}
\author{A.\,V.\ Vasil'ev}
\address{Sobolev Institute of Mathematics, Novosibirsk, Russia;\newline
and Novosibirsk State University, Novosibirsk, Russia}
\email{vasand@math.nsc.ru}
\author{E.\ Vdovin}
\address{Sobolev Institute of Mathematics, Novosibirsk, Russia;\newline
and Novosibirsk State University, Novosibirsk, Russia}
\email{	vdovin@math.nsc.ru}

\thanks{O'Brien was supported by the Marsden Fund of New Zealand
grant UOA 107; Ponomarenko and Vasil'ev were supported
by the Mathematical Center in Akademgorodok under agreement
No.\ 075-15-2019-1613 with the Ministry of Science and Higher
Education of the Russian Federation; Vdovin was supported by the 
RFBR grant No.\ 18-01-00752. We thank the referee for 
helpful feedback.}
\date{}

\begin{abstract}
Let $m$ be a positive integer and let $\Omega$ be a finite set.
The $m$-closure of $G\le\sym(\Omega)$ is the largest
permutation group on $\Omega$ having the same orbits as $G$ in its
induced action on the Cartesian product~$\Omega^m$.
The $1$-closure and $2$-closure of a solvable permutation group
need not be solvable. We prove that the $m$-closure of a solvable permutation
group is always solvable for $m\ge 3$.
\end{abstract}

\maketitle

\centerline{\it Dedicated to the memory of our friend Jan Saxl}

\section{Introduction}
Let $m$ be a positive integer and let $\Omega$ be a finite set.
The {\it $m$-closure} $G\m$ of $G\le\sym(\Omega)$ is
the largest permutation group on $\Omega$ having
the same orbits as $G$ in
its induced action on the Cartesian product~$\Omega^m$.
Wielandt \cite[Theorems~5.8 and 5.12]{Wielandt1969} showed that
\qtnl{100120p}
G^{(1)}\ge G^{(2)}\ge\cdots\ge G^{(m)}=G^{(m+1)}=\cdots =G,
\eqtn
for some $m<|\Omega|$.
(Since the stabilizer in $G$ of all but one point is always trivial,
$G^{(n-1)}=G$ where $n = | \Omega |$; see Theorem~\ref{281219a}.)
In this sense, the $m$-closure can be considered as a natural
approximation of~$G$.
Here we study the closures of
solvable groups; for the nonsolvable case,
see~\cite{LPS1988,PS1992,XuGLP2011}.

The $1$-closure of $G$ is the direct product of symmetric
groups $\sym(\Delta)$, where~$\Delta$ runs over the orbits of~$G$.
Thus the $1$-closure of a solvable group is solvable if and only
if each of its orbits has cardinality at most~$4$.
The case of $2$-closure is more interesting.
The $2$-closure of every (solvable) $2$-transitive group
$G\le\sym(\Omega)$ is $\sym(\Omega)$; other examples of
solvable~$G$ and nonsolvable $G^{(2)}$ appear in~\cite{Skresanov2019}.
But, as shown by Wielandt ~\cite{Wielandt1969},
each of the classes of finite $p$-groups and groups of odd order
is closed with respect to taking the $2$-closure.
Currently, no characterization of solvable groups having
solvable  $2$-closure is known.

Seress~\cite{Sere1996} observed that if $G$ is a primitive solvable group, then
$G^{(5)}=G$; so the $5$-closure of a primitive solvable group is
solvable. Our main result is the following stronger statement.

\thrml{241219a}
The $3$-closure of a solvable permutation group is solvable.
\ethrm

Theorem~\ref{241219a} follows from
Theorems~\ref{220520d}, \ref{110519c}, and~\ref{301219b}.
The corollary below  is an  immediate consequence of
Theorem~\ref{241219a} and the chain of inclusions~\eqref{100120p}.

\crllrl{241219b}
For every integer $m\ge 3$, the $m$-closure of a solvable
permutation group is solvable.
\ecrllr

We briefly outline the structure of our proof.
In Section~\ref{090120c} we recall the
basic theory of the closure of permutation groups,
as developed by Wielandt~\cite{Wielandt1969}.
In Section~\ref{090120b} we
deduce that the $m$-closure of the
direct (respectively, imprimitive wreath) product of two permutation
groups is isomorphic to a subgroup of the direct
(respectively, imprimitive wreath) product of their $m$-closures,
and prove (with some natural constraints) that the same holds true for
the primitive wreath product.
Thus, in Theorem~\ref{220520d}, we reduce the proof of
Theorem~\ref{241219a} from an arbitrary solvable permutation group $G$
to a linearly primitive group: a point stabilizer $G_0$ of
such a group is a primitive linear group over a finite field.

A natural dichotomy arises in our treatment of linearly primitive groups.
If the point stabilizer $G_0$ has a regular faithful orbit,
then Wielandt's theory shows that the $3$-closure of $G$ is solvable (Corollary~\ref{060519b}).
Otherwise, we use results from~\cite{Yang2010,Yang2011,Yang2016} to obtain
in Theorem~\ref{110519c}
an explicit list of pairs $(d,p)$ for which
$G\le\AGL(d,p)$.

In Section~\ref{090120h} we complete the proof, relying on computer
calculations.
Namely, following Short's approach~\cite{Short1992},
we find for each pair $(d,p)$ a set $\cH$ of linearly primitive subgroups of $\GL(d,p)$ containing
a conjugate of each maximal solvable primitive subgroup  of $\GL(d,p);$ in particular, $G_0$ is a subgroup
of some $H\in\cH$. If $G$ is $2$-transitive, then $G^{(3)}$ is solvable (Lemma~\ref{220520b}). In the remaining cases, $G^{(3)}=G$ (and so solvable) because $H^{(2)}=H$ for every $H\in\cH$. Verification of the latter is based on sufficient conditions given in Corollary~\ref{060519b} and Lemma~\ref{220520c}.

\section{Wielandt's theory}\label{090120c}
Let $G\le\sym(\Omega)$ and let $m$ be a positive integer. The $m$-orbits of
$G$ are the orbits of componentwise action of~$G$ on the Cartesian
product~$\Omega^m$ of~$\Omega$; the set of all such orbits is denoted
by $\orb_m(G)$.

\xmpll{100620a}
We describe the $m$-orbits of $G=\sym(\Omega)$.
Given an $m$-tuple $\alpha=(\alpha_1,\ldots,\alpha_m)$ of~$\Omega^m$,
let $\pi(\alpha)$ be the partition of
$I=\{1,\ldots,m\}$ such that the elements $i$ and $j$ belong to the
same class of~$\pi$ if and only if $\alpha_i=\alpha_j$.
If $s\in\orb_m(G)$, then 
$\pi(s) :=\pi(\alpha)$ does not depend on the choice of $\alpha\in s$.
If $m\le|\Omega|$, then the mapping $s\mapsto \pi(s)$ establishes a
one-to-one correspondence between the $m$-orbits of $G$ and the
partitions of~$I;$ in particular, $|\orb_m(G)|\ge m+2$ for all $m\ge 3$.
\exmpl

A permutation group $H$ on $\Omega$ is {\it $m$-equivalent} to~$G$ if
$$
\orb_m(H)=\orb_m(G).
$$
Obviously, $m$-equivalence is an 
equivalence relation on the set of permutation groups on~$\Omega$.
If $m\ge 2$, then two $m$-equivalent groups share some properties,
such as primitivity, or $2$-transitivity;
see \cite[Theorems~4.8 and~4.10]{Wielandt1969}.

The following criterion for $m$-equivalence can
easily be deduced from \cite[Theorem~4.7]{Wielandt1969}.
\lmml{07052020a}
Let $G$ and $H$ be permutation groups  on $\Omega$ and assume that $G\le H$.
Then $H$ is $m$-equivalent
to $G$ if and only if, for every $\alpha\in\Omega^m$ and every $h\in H$,
there exists $g\in G$ such that $\alpha^h=\alpha^g$.
\elmm

Wielandt \cite[Theorem~4.3 and Lem\-ma~4.12]{Wielandt1969} established
the following.

\thrml{100120a}
Let $m\geq2$ be an integer and let $G$ and $H$ be $m$-equivalent permutation groups
on~$\Omega$. The following hold:
\nmrt
\tm{i}  $G$ and $H$ are $(m-1)$-equivalent;
\tm{ii} $G_\alpha$ and $H_\alpha$ are $(m-1)$-equivalent for all $\alpha\in\Omega$.
\enmrt
\ethrm

The definition of $m$-equivalence implies that the
{\em $m$-closure} of $G$ is the largest group in
the class of $m$-equivalent groups containing~$G$. In particular,
$G$ and $H$ are $m$-equivalent if and only if $G\m=H\m$.
If $G\m=G$ then $G$ is {\it $m$-closed}.
Note that $m$-closure is a
closure
operator: namely,
$G\le G\m$, $G\m=(G\m)\m$, and $G\le H$ implies $G\m\le H\m$.

\thrml{281219a}
Let $m\ge 2$ be an integer.
 If a point stabilizer of a permutation group $G$ is $(m-1)$-closed,
then $G$ is $m$-closed.
\ethrm
\proof Let $G\le\sym(\Omega)$  and let $\alpha\in\Omega$.
Since $G$ and $H:=G\m$ are $m$-equivalent, $G_\alpha$
and $H_\alpha$ are $(m-1)$-equivalent (Theorem~\ref{100120a}(ii)).
Since $G_\alpha$ is $(m-1)$-closed,
$$
H_\alpha\le (H_\alpha)^{(m-1)}=(G_\alpha)^{(m-1)}=G_\alpha.
$$
But
$G_\alpha\le H_\alpha$, so $G_\alpha=H_\alpha$. Furthermore,
Theorem~\ref{100120a}(i) implies that $G$ and $H$ are $1$-equivalent.
Therefore $\alpha^G=\alpha^H$.
Hence
$$
|G|=|\alpha^G|\cdot|G_\alpha|=|\alpha^H|\cdot|H_\alpha|=|H|.
$$
Thus $G=G\m$ because $G\le H$.
\eprf\medskip

A permutation group is {\it partly regular} if it has a faithful regular orbit. Clearly, every subgroup of a partly regular group is partly regular.

\crllrl{060519b}
Let $G\le\sym(\Omega)$ and let $m\ge 2$ be an integer.
If an $(m-1)$-point stabilizer of $G$ is partly regular, then $G$
is $(m+1)$-closed.
\ecrllr
\proof If $G$ has a partly regular $(m-1)$-point stabilizer, then $G$ has
an $m$-point stabilizer which is trivial and so $m$-closed.
The claim follows from Theorem~\ref{281219a}.
\eprf

%
%

\section{Closures of permutation groups}\label{090120b}
We now study the $m$-closure operator under
standard operations in permutation group theory;
compare with similar results for $m=2$ in~\cite{EP01}.

\thrml{241219c}
Let $K\le\sym(\Gamma)$,  let $L\le\sym(\Delta)$, and  let $K\times L$ act
on the disjoint union $\Gamma\cup\Delta$. For every integer $m\ge 1$,
$$
(K\times L)\m\le K\m \times L\m.
$$
\ethrm
\proof
Observe that $K\times L$ and $H:=(K\times L)\m$ are
$1$-equivalent. Hence the sets $\Gamma$ and $\Delta$ are
invariant under~$H$. It follows that
$K=(K\times L)^{\Gamma}$ is $m$-equivalent to~$H^{\Gamma}$,
the permutation group induced by the action of $H$ on $\Gamma$,
 and
$L=(K\times L)^{\Delta}$ is $m$-equivalent to $H^{\Delta}$. 
In particular, $H^\Gamma\le K\m$ and $H^\Delta\le L\m$. Thus
$$
(K\times L)\m=H\le H^{\Gamma}\times H^{\Delta}\le  K\m \times L\m,
$$
as required.
\eprf\medskip

An analogue of Theorem~\ref{241219c} exists for the direct
product of permutation groups acting on the Cartesian product of
their underlying sets, but it is not needed here.

The following is a consequence of~\cite[Lemma~2.5]{KaluK1976}.
\thrml{241219d}
Let $K\wr L$ be the imprimitive wreath product of permutation
groups $K$ and $L$.  For every integer $m\ge 2$,
$$
(K\wr L)\m=K\m \wr L\m.
$$
\ethrm

The case of the wreath product in {\it product action} is more subtle.
Recall, for example from \cite[Lemma~2.7A]{DM}, that
$K\uparrow L$,
the wreath product in product action of permutation groups $K$ and $L$,
is primitive
if and only if $K$ is primitive and nonregular, and $L$ is transitive and
nontrivial.
For the remainder of the paper,
we assume that
$K\uparrow L$ {\it is primitive} and so label this construction as
the {\it primitive wreath product}.
Even with this assumption, $(K\uparrow L)\m$ is not
always a subgroup of
$K\m \uparrow L\m$: consider for example $m=2$, $K=\sym(4)$, and $L=\alt(3)$.
But we obtain the following.

\thrml{241219d1}
Let $K\uparrow L$ be the primitive wreath product of permutation
groups~$K$ and $L$ of degrees $r$ and $d$, respectively.
Assume that $m\ge 3$ is an integer such that $m\le r$,
and also $m\le d$ unless $d=2$. Then
$$
(K\uparrow L)\m\le K\m \uparrow L\m.
$$
\ethrm

\proof Let $K\le\sym(\Gamma)$ and $L\le\sym(\Delta)$, where
$\Gamma$ and $\Delta$ are sets of cardinality $r$ and~$d$, respectively.
Without loss of generality, we assume that  $\Delta=\{1,\ldots,d\}$.
Thus $G:=K\uparrow L$ acts on the Cartesian product
$$
\Omega=\underbrace{\Gamma\times \cdots\times \Gamma}_{d \ \,\text{copies}}.
$$
In what follows, an $m$-tuple $x\in\Omega^m$ is considered as an $m\times d$ matrix $(x_{ij})$ with $x_{ij}\in \Gamma$; the $j$th column of $x$ is denoted by~$x_{*j}$. Note that $K^d$ acts on $x$ by permuting elements inside columns, whereas $L$ permutes the columns.

As observed in \cite[Proof of Proposition~3.1]{EP01},
$G$ is contained in the $2$-closed group
$\sym(\Gamma)\uparrow\sym(\Delta)$. It follows that $H:=G\m$,
being $2$-equivalent to~$G$, is also contained in
$\sym(\Gamma)\uparrow\sym(\Delta)$. Therefore every
permutation  of~$H$ can be written in the form
\qtnl{070120a}
h=(h_1,\ldots,h_d;\ov h),\qquad h_1,\ldots,h_d\in\sym(\Gamma),\ \ov h\in\sym(\Delta).
\eqtn
Let $\ov H = \{ \ov h \;:\; h\in H\}$ be the permutation group
induced by the action of~$H$ on~$\Delta$.

\noindent
As a critical step in our proof, we establish the following.

\medskip
\noindent
{\bf Claim.} {\it $H\le \sym(\Gamma)\uparrow L\m$.}\medskip

\proof  Without loss of generality, we may assume that $d\ge 3$.
It suffices to show that $\ov H\le L\m$.
Equivalently (cf.\ Lemma~\ref{07052020a}), we show that, for every
$\alpha \in\Delta^m$ and every $\ov h\in \ov H$,
there exists $\ov g\in L$ such that
\qtnl{200220a}
\alpha^{\ov h}=\alpha^{\ov g}.
\eqtn
Since by assumption $K \uparrow L$ is primitive,
$L$ is transitive
and a subgroup of $\ov{H}$. Consequently, $\ov H=L\ov{H}_d$, where
$\ov{H}_d$ is the stabilizer in~$\ov H$ of the point $d\in \Delta$.
Thus we may assume that the element $\ov h$ in
\eqref{200220a} belongs to $\ov{H}_d$.

Let $\alpha=(\alpha_1,\ldots,\alpha_m)\in\Delta^m$, 
where $\{\alpha_1,\ldots,\alpha_m\}=\{j_1,\ldots,j_t\}$ 
and $1\le j_1<\cdots<j_t\le m$.
Let $\pi(\alpha)$ be  the partition of 
$\{1,\ldots,m\}$ into the classes $\{\ell:\ \alpha_\ell=j_i\}$ where 
$i=1,\ldots,t$ and $t=|\pi(\alpha)|$.
Note that the number of $m$-orbits of $\sym(\Gamma)$ is at least $m+2$, because $|\Delta|=d\ge m \ge 3$ (see Example \ref{100620a}). Thus there are  pairwise distinct $m$-orbits $s_0,s_1,\ldots,s_t$  such that $\pi(s_i)\ne\pi(\alpha)$ for all~$i$.
Denote by $\myx(\alpha)$ the set of all tuples $x\in\Omega^m$ such that
\qtnl{120621a}
\pi(x_{*j_1})=\pi(s_1),\quad \pi(x_{*j_2})=\pi(s_2),
\quad\ldots,\quad
\pi(x_{*j_t})=\pi(s_t),
\eqtn
\qtnl{120621b}
\pi(x_{*d})=\pi(\alpha),\quad\pi(x_{*i})=\pi(s_0),\quad i\not\in\{j_1,\ldots,j_t,d\}.
\eqtn
Of course, the choice of $s_0,\ldots,s_t$ is arbitrary. 
However, if $\beta$ is in the $m$-orbit of $\sym(\Delta)$ containing $\alpha$, then $\pi(\beta)=\pi(\alpha)$, so we can use the same $s_0,\ldots,s_t$ to define $\myx(\beta)$.

Assume that $\myx(\alpha)=\myx(\beta)$ for some $\beta\in\Delta^m$.
Condition~\eqref{120621b} ensures that
$$
\pi(\alpha)=\pi(\beta)\qaq \{\alpha_1,\ldots,\alpha_{m-1}\}=
\{\beta_1,\ldots,\beta_{m-1}\}.
$$
In particular, $\alpha_i=\alpha_j$ if and only
if $\beta_i=\beta_j$ for all $1\le i,j\le m$.
Condition~\eqref{120621a} yields $\alpha=\beta$.
Thus
\qtnl{270520d}
\myx(\alpha)\cap \myx(\beta)\neq\varnothing\quad
\text{if and only if}\quad \alpha=\beta.
\eqtn

Let $h\in H$ be as in ~\eqref{070120a} and let $\ov{h}\in\ov{H}_d$.
In view of our definitions, given
indices $i,j\in\{1,\ldots,m\}$, the column $j$ of a
tuple $x\in \myx(\alpha)$ belongs to~$s_i$ if and only if the
column $j\phmb{\ov h}$ of the tuple $x^h$ belongs to~$s_i$. Therefore
\qtnl{190220a}
\myx(\alpha)^h=\myx(\alpha^{\ov h}).
\eqtn
In particular, the preimage of $\ov{H}_d$ in $H$ acts on the
set $\myx(\alpha)$ for each $\alpha\in\Delta^m$. Since $G\le H$,
\eqref{190220a} holds for every $g\in G$
with $\ov{g}\in L_d$, where $L_d$ is the stabilizer in $L$ of the
point $d\in\Delta$.

Since $G$ and $H$ are $m$-equivalent and $G\le H$,
Lemma~\ref{07052020a} implies that for
every $x\in \myx(\alpha)$ and every $h\in H$ there
exists $g\in G$ such that $x^h=x^g$.
But $\ov h\in\ov{H}_d$; so, in accord
with the first equality in~\eqref{120621b}, $\ov{g}\in L_{d}$. In view
of \eqref{190220a}, this implies that
$$
x^g=x^h\in \myx(\alpha^{\ov h})\cap \myx(\alpha^{\ov g}).
$$
Now \eqref{270520d} yields $\alpha^{\ov{h}}=\alpha^{\ov{g}}$,
which proves \eqref{200220a}. \eprf\medskip

We return to the proof of the theorem.
It remains to verify that for every $h\in H$ and every $\alpha\in\Gamma^m$,
\qtnl{270520e}
\alpha^{h_j}\in \alpha^{K}\mbox{ for }j=1,\ldots,d,
\eqtn
where $h_j$ is defined by \eqref{070120a}.
Without loss of generality, we prove~\eqref{270520e} for $j=1$ only.
Take distinct $m$-orbits $s_0$ and $s_1$ of $\sym(\Gamma)$ such
that $\alpha_1\in s_1$. Choose a tuple $x\in\Omega^m$ satisfying
the conditions
\qtnl{060620i}
x_{*1}=\alpha_1\qaq x_{*2},\ldots,x_{*d}\in s_0.
\eqtn

Our claim implies
that $L\le \ov H\le L\m$, so
$L$ and $\ov H$ are $m$-equivalent and  have the
same orbits. Hence there exists $\ov g\in L$ such that
$\ov h\ov g\in\sym(\Delta)$ fixes $1\in\Delta$. So
$$
g=(1,\ldots,1;\ov g) \in G
$$
and the only column of the tuple~$x^{hg}$ belonging to $s_1$ is the
first one, from \eqref{060620i}. Thus
$$
(x^{hg})_{*1}=\alpha^{h_1}.
$$
Since $x^H\in\orb_m(H)=\orb_m(G)$, there exists $f\in G$ such
that $x^{hg}=x^f$; in particular, $\ov f\in\sym(\Delta)$ fixes $1\in\Delta$.
Consequently,
$$
\alpha^{h_1}=(x^{hg})_{*1}=(x^f)_{*1}=\alpha^{f_1}\in \alpha^{K},
$$
as required.\eprf\medskip

If $m=3$, then the hypothesis of Theorem~\ref{241219d1} does not impose any
restrictions on the degrees $r$ and $d$ of the groups $K$ and $L$ because
the primitivity of
$K\uparrow L$ implies that
$r\ge3$ and $d\ge2$.
Therefore the following holds.

\crllrl{211220a}
Let $K\uparrow L$ be the primitive wreath product of permutation
groups~$K$ and $L$. Then
$$
(K\uparrow L)\3\le K\3 \uparrow L\3.
$$
\ecrllr

\section{Reducing primitive solvable groups to
linearly primitive groups}\label{280520b}

We summarize the well-known structure of a
primitive solvable group; for a proof, see
for example \cite[Chap.\ 1, Theorem~7]{Su76}.

\thrml{140311c}
Let $G\le\sym(\Omega)$ be a primitive solvable permutation group.
Now $\Omega$ has cardinality $p^d$
for a prime $p$ and integer $d\ge 1$ and
can be identified with a $d$-dimensional vector space $V$ over $\GF(p)$.
Moreover,  $G\le\AGL(d,p)$,
and the stabilizer $H$ in $G$ of the zero vector is an irreducible
subgroup of $\GL(d,p)$.
\ethrm

We establish Theorem \ref{241219a} for two classes of groups.
\lmml{220520a}
If $G\leq\AGaL(1,p^d)$, then $G^{(3)}$ is solvable.
\elmm

\proof
A point stabilizer $\GaL(1,p^d)$
of $\AGaL(1,p^d)$ 
is $2$-closed by \cite[Proposition 3.1.1]{XuGLP2011}.
Theorem~\ref{281219a} implies that $\AGaL(1,p^d)$ is $3$-closed.
Thus
$$
G^{(3)}\le \AGaL(1,p^d)^{(3)}=\AGaL(1,p^d),
$$
and so $G^{(3)}$ is solvable.
\eprf

\lmml{220520b}
If $G$ is a $2$-transitive solvable group, then $G^{(3)}$ is solvable.
\elmm

\proof
By a theorem of Huppert \cite[Theorem~6.9]{ManzWolf1993},
$G\le\AGaL(1,p^d)$, or
\qtnl{311219a}
p^d\in\{3^2,5^2,7^2, 11^2,23^2,3^4\}.
\eqtn
The first case is settled by Lemma~\ref{220520a}.

In the second case, we consider a point stabilizer $G_\alpha$ and used the {\tt TwoClosure} command
in {\sf GAP} \cite{gap}
and the IRREDSOL package \cite{Hoef2017}
to establish that $(G_\alpha)^{(2)}$ is solvable for
those $p^d$ satisfying~\eqref{311219a}.
Since $(G^{(3)})_\alpha$ and $G_\alpha$ are $2$-equivalent,
$$
(G^{(3)})_\alpha\le (G_\alpha)^{(2)},
$$
and so $(G^{(3)})_\alpha$ is solvable.
Hence $G^{(3)}$, an
extension of an elementary abelian group (of order $p^d$) by
$(G^{(3)})_\alpha$, is solvable.\eprf\medskip

The following lemma gives a sufficient condition for
a primitive solvable group to be $3$-closed.

\lmml{220520c}
In the notation of Theorem~$\ref{140311c}$,
if there exists a nonzero $\alpha \in V$ such that the restriction of
$H$ to $\alpha^H$ is $2$-closed, then $G$ is $3$-closed.
\elmm

\proof
Since $H$ is an irreducible linear group,
the orbit $\Delta=\alpha^H$ contains a basis of $V$.
It follows that $H\cong H^\Delta$,
the restriction of $H$ to $\Delta$.
On the other hand, $G^{(3)}$ and $G$ are $2$-equivalent
by Theorem~\ref{100120a}(i). Consequently, $G^{(3)}$ is also a
primitive subgroup of $\AGL(d,p)$ \cite[Theorem~1.4]{XuGLP2011}.
Hence the stabilizer $L$ in $G^{(3)}$ of the zero vector acts
faithfully on~$\Delta$. Thus
\qtnl{070620a}
H\cong H^\Delta\qaq L\cong L^\Delta.
\eqtn
Since $H$ and $L$ are point stabilizers of $3$-equivalent groups, they
are $2$-equivalent by Theorem~\ref{100120a}(ii). Therefore
\qtnl{070620b}
L\le H^{(2)}.
\eqtn
Following \cite[Lemma~2.1(iii)]{PonomarenkoV2020}, we verify that
\qtnl{070620c}
(H^{(2)})^\Delta\le (H^\Delta)^{(2)}.
\eqtn
By hypothesis $(H^\Delta)^{(2)}=H^\Delta$.
This,
together with \eqref{070620a}, \eqref{070620b}, and \eqref{070620c},
implies that
$$
L\cong  L^\Delta\le (H^{(2)})^\Delta\le (H^\Delta)^{(2)}=H^\Delta\cong H.
$$
Since $H\le L$, this yields $H=L$. Thus $G=G^{(3)}$.
\eprf
\medskip

If the point stabilizer $H$ in Theorem~\ref{140311c} is primitive
as a linear group, then $G$
is {\it linearly primitive}, otherwise it is {\it linearly imprimitive}.
As is well-known, the latter case reduces to the primitive wreath product;
see for example~\cite[Proposition~4.1]{EP01}.

\lmml{160311a}
Every linearly imprimitive solvable permutation group $G$ is
isomorphic to a subgroup of a primitive wreath product of two
solvable permutation groups of degrees smaller than the degree of~$G$.
\elmm

We now reduce the proof of Theorem \ref{241219a} to linearly primitive groups.

\thrml{220520d}
A counterexample of minimal degree to the statement
of Theorem~$\ref{241219a}$ is linearly primitive.
\ethrm

\proof
We consider separately the cases where the permutation group $G$ is
intransitive, imprimitive, or linearly imprimitive.\medskip

{\bf Case 1:} $G$ is intransitive. Now $G$ is a subdirect product
of (solvable) constituents, say $K$ and $L$, and their degrees are
less than that of $G$. Therefore their $3$-closures are solvable.
By Theorem~\ref{241219c} so is $(K\times L)^{(3)}$. Thus
$$
G^{(3)}\le (K\times L)^{(3)}
$$
is also solvable.\medskip

{\bf Case 2:} $G$ is imprimitive. Now $G$ can be identified with a
subgroup of the imprimitive wreath product $K\wr L$, where $K$ and $L$
are solvable permutation groups, and their degrees are less than
that of~$G$. Therefore their $3$-closures are solvable. By
Theorem~\ref{241219d} so is $(K\wr L)^{(3)}$. Thus
$$
G^{(3)}\le (K\wr L)^{(3)}
$$
is also solvable.\medskip

{\bf Case 3:} $G$ is primitive, but linearly imprimitive.
Now, by Lemma~\ref{160311a}, it can be identified with a subgroup of
the primitive wreath product $K\uparrow L$, where $K$ and~$L$ are
solvable permutation groups, and their degrees are less than
that of~$G$. Therefore their $3$-closures are solvable.
By Corollary~\ref{211220a} so is $(K\uparrow L)^{(3)}$.
Thus
$$
G^{(3)}\le (K\uparrow L)^{(3)}
$$
is also solvable.
\eprf

\section{Linearly primitive solvable groups: background theory}\label{280520c}

Let $G$ be a linearly primitive solvable permutation group with
point stabilizer $G_0$.  Often essential information about $G_0$
can be obtained from a maximal solvable primitive
linear group $H$ containing $G_0$.
Basic information on the structure of $H$
is collected in the following theorem; for its proof, see for example
\cite[Lemma~2.2]{Sere1996}.

\thrml{060519a}
Let $H\le\GL(d,p)$ be a maximal solvable primitive group.
It has a  series $1<U\le F\le A\le H$  satisfying the following:
\begin{enumerate}
\item[(i)]
$U$ is the unique maximal abelian normal subgroup of $H$, the linear span of $U$ in {\em Mat}$(d,p)$ is $\GF(p^a)$ where $a$ divides~$d$, and $U$ is cyclic of order $p^a-1;$
\item[(ii)]
$F=\Fit(C_H(U))$ is the Fitting subgroup of the centralizer $C_H(U)$, and
\mbox{$|F\,/\,U|$}$\,=\,e^2$ where $d = ae$ and
each prime divisor of $e$ divides $p^a-1;$
\item[(iii)]
$A=C_H(U)$ and $A/F$ is isomorphic to a completely reducible subgroup
of the direct product $\prod_{i=1}^m\Sp(2n_i,p_i)$ where the $p_i$
and $n_i$ are defined by the
prime power decomposition $e=\prod_{i=1}^mp_i^{n_i};$
\item [(iv)]
$H/A$ is isomorphic to a subgroup of $\aut(\GF(p^a))$ and so $|H/A|$ divides $a$.
\end{enumerate}
\ethrm

Let $G$ be a linearly primitive solvable permutation group.
We fix an embedding of
a point stabilizer $G_0$ of $G$ into a specific
maximal solvable primitive linear group $H$.
We call the integers $p$, $d$, $a$, $e$
defined in Theorem~\ref{060519a} for $H$
the {\it parameters} of~$G$.
Although the parameters of $G$ depend on the choice of $H$,
the monotonicity of the $m$-closure operator guarantees that our
results are independent of it.

\thrml{110519c}
The $3$-closure of a linearly primitive solvable permutation group is solvable,
except possibly for those groups whose parameters are
listed in columns~$2$--$5$ of Table~$\ref{t2}$.
\ethrm
\proof Let  $G$ be a linearly primitive solvable permutation group with
parameters $p$, $d$, $a$, $e$, and let $H\le\GL(d,p)$ be
a maximal solvable primitive group containing a point
stabilizer $G_0$ of~$G$.

If $e=1$, then $G\le\AGL(1,p^d)$, and
the result follows
by Lemma~\ref{220520a}. So we may assume that $e>1$.
If $H$ is partly regular, then so is $G_0\le H$;
thus $G$ is $3$-closed by Corollary~\ref{060519b} (take $m=2$).
If $H$ is not partly regular,  then its parameters
are listed in \cite[Corollary~3.2]{Yang2016};
columns~$2$--$5$ of Table~$\ref{t2}$ are
taken from~\cite[Table~2]{Yang2016}.~\eprf\medskip

The data in \cite[Table~2]{Yang2016}
was obtained using
\cite[Theorem~4.1]{Yang2011} which states the following: if $H$
is not partly regular, then
\qtnl{301219a}
e=2,\ 3,\ 4,\ 8,\ 9,\ 16.
\eqtn
Hence $e = r^k$ where $2\leq r \leq 3$ and $1\leq k \leq 4$.
Let $b$ be the least positive integer with $p^b\equiv1\pmod{r^c}$ where $c=2$ for $r=2$ and $c=1$ otherwise. We denote the general orthogonal group of degree $2k$ over the field of order $q$ by $O^\varepsilon(2k,q)$ where $\varepsilon\in\{+,-\}$ depends on the Witt index of the corresponding quadratic form. The {\em $r$-radical} of a group is its largest normal $r$-subgroup.

\lmml{211220b} We retain the notation of Theorem~$\ref{060519a}$ where $e=r^k$.
\begin{enumerate}
\item[(i)]
$F$ is the central product of $U$ and an extraspecial group $E$
of order $r^{2k+1};$ if $b \,|\, a$, then all such subgroups
$F$ are conjugate in $\GL(d,p)$,
else $r=b=2$, $a$ is odd, and there are two conjugacy classes of such subgroups;
\item[(ii)]
$A/F$ is a maximal solvable subgroup of $N/F$, where $N=N_L(F)$ and $L=C_{\GL(d,p)}(U)\cong\GL(e,p^a);$ if $b \,|\, a$, then $N/F\cong\Sp(2k,r)$,
else $N/F$ is isomorphic to one of $O^\varepsilon(2k,2)$ depending on the conjugacy class of $F;$
\item[(iii)]
the $r$-radical of $A/F$ has order at most $2$ and is trivial if $b \,|\, a$.
\end{enumerate}
\elmm

\proof Item (i) is well known; see for example
\cite[Theorem~2.4.7]{Short1992}
and subsequent remarks. 
By \cite[Theorems~2.5.31, 2.5.34 and~2.4.12]{Short1992},
$A/F$ is isomorphic to a maximal solvable subgroup $M$ of $S=N/F$,
where $S=\Sp(2k,r)$ if $b\,|\,a$ and $S$ is one of
$O^\varepsilon(2k,2)$ otherwise; this proves (ii).
Moreover, $M$ fixes no nonzero isotropic subspace of the
natural $\GF(r)$-module of~$S$.
Therefore $S$ is not contained in any parabolic subgroup of $S$,
so the $r$-radical of $M$ is either trivial or has order $2$,
and the latter is possible only if $M$ is an orthogonal group.~\eprf\medskip

The following lemma can be established using {\sc
Magma} \cite{magma}. Note that $O^+(2,2)$,
$\Sp(2,2)\cong O^-(2,2)$, $\Sp(2,3)$,  and $O^+(4,3)$ are solvable.

\lmml{spo}
Let $M$ be a maximal solvable subgroup of
$S \leq \GL(d, r)$ and let the $r$-radical of $M$ satisfy the
order conditions of Lemma~{\em\ref{211220b}(iii)}.
\nmrt
\item[(i)] If $S=\Sp(2,2),\,\Sp(2,3),\,O^\varepsilon(2,2),\,O^+(4,3)$,
then $M=S$.

\item[(ii)] If $S=O^-(4,2)$, then $M$ is conjugate to
a subgroup isomorphic to either $5:4$ of order $20$,
or $S_3\times S_2$  of order $12$.

\item[(iii)] If $S=\Sp(4,2)$, then $M$ is conjugate to
a subgroup isomorphic to either
$O^+(4,2)\cong S_3\wr S_2$ of order $72$, or the normalizer
of a Sylow $5$-subgroup
of $\Sp(4,2)$ which is isomorphic to $5:4$ and has order $20$.

\item[(iv)] If $S=\Sp(4,3)$, then $M$ is conjugate to
one of the following:
\begin{itemize}
\item the normalizer of a Sylow $5$-subgroup of $\Sp(4,3)$ which
is isomorphic to $D_{20}.2$ and has order $40$;
\item a subgroup isomorphic to $2^{1+4}:S_3$ of order $192$;
\item a subgroup isomorphic to $2^{1+4}:D_{10}$ of order $320$;
\item a subgroup isomorphic to $\Sp(2,3)\wr S_2$ of order $1152$.
\end{itemize}

\item[(v)] If $S=O^+(6,2)$, then $M$ is conjugate to one of the following:
\begin{itemize}
\item
the normalizer of a Sylow $3$-subgroup of $O^+(6,2)$
which is isomorphic to
$O^+(4,2)\times O^+(2,2)\cong (S_3\wr S_2)\times S_2$
and has order $144$;
\item
the normalizer of a Sylow $5$-subgroup of $O^+(6,2)$ which
is isomorphic to $(5:4)\times S_3$ and has order $120$;
\item
the normalizer of a Sylow $7$-subgroup of $O^+(6,2)$
which is isomorphic to $7:6$ and has order $42$.
\end{itemize}

\item[(vi)] If $S=O^-(6,2)$, then $M$ is conjugate to one of the following:
\begin{itemize}
\item a subgroup isomorphic to $3^{1+2}:(2.S_4)$ of order $1296$;
\item a subgroup  isomorphic to $3^3:(S_4\times S_2)$  of order $1296$;
\item the normalizer of a Sylow $5$-subgroup of $O^-(6,2)$ which is
isomorphic to $(5:4)\times S_2$ and has order $40$.
\end{itemize}

\item[(vii)] If $S=\Sp(6,2)$, then $M$ is conjugate to
one of the following:
\begin{itemize}
\item a subgroup isomorphic to $3^{1+2}:(2.S_4)$ of order $1296$;
\item a subgroup  isomorphic to $3^3:(S_4\times S_2)$  of order $1296$;
\item the normalizer of a Sylow $5$-subgroup of $\Sp(6,2)$ which
is isomorphic to $(5:4)\times S_3$ and has order $120$;
\item the normalizer of a Sylow $7$-subgroup of $\Sp(6,2)$ which
is isomorphic to $7:6$ and has order $42$.
\end{itemize}
\enmrt
\elmm


In Table~\ref{t1},
we summarize the orders of the maximal solvable subgroups listed in
Lemma \ref{spo}. 

\begin{table}[h]
\begin{center}
\caption{The orders of certain maximal solvable subgroups}\label{t1}
\begin{tabular}{|c|c|l|}
\hline
$e$        & $S$ & $|M|$  \\
\hline
$9$        & $\Sp(4,3)$            & $40,  192, 320,  1152$ \\
\hline
$8$       & $\Sp(6,2)$             &  $42,120,1296$ \\
          & $O^+(6,2)$             & $42, 120, 144$ \\
          & $O^-(6,2)$             & $40, 1296$ \\ \hline
$4$       & $\Sp(4,2)$             &  $20,72$ \\
          & $O^+(4,2)$             & $72$ \\
          & $O^-(4,2)$             & $12,20$\\
\hline
$3$       & $\Sp(2,3)$             & $24$ \\
\hline
$2$       & $\Sp(2,2)$             & $6$ \\
          & $O^+(2,2)$             & $2$ \\
          & $O^-(2,2)$             & $6$ \\
\hline
\end{tabular}
\end{center}
\end{table}

To state the next lemma we introduce some additional notation.
Let $I_m$ be the identity $m\times m$ matrix, and let $\otimes$
denote the Kronecker product of matrices. If $y$ and $z$ are
$k\times k$ and $m\times m$ matrices respectively,
then $y\otimes I_m$ commutes with $I_k\otimes z$. 
For $G\leq \GL(k,p)$ let $G\otimes I_m$ be the
subgroup $\{g\otimes I_m\mid g\in G\}$ of $\GL(k\cdot m,p)$.
Clearly $G\cong G\otimes I_m$.
For
positive integers $k,m$ there exists a natural embedding
\qtnl{GLextfield}
\GL(k,p^m):\langle \varphi\rangle \leq \GL(k\cdot m,p),
\eqtn
where $\varphi$ is a field automorphism of $\GL(k,p^m)$.
The uniqueness of a field of order $p^m$ implies that all
such embeddings are conjugate in $\GL(k\cdot m,p)$. Below we
assume that we fix such an embedding and so realise
$\GL(k,p^m):\langle \varphi\rangle$ as a subgroup of $\GL(k\cdot m,p)$.

\lmml{Kron} We retain the notation of Theorem~$\ref{060519a}$ where $e=r^k$.
\nmrt
\item[(i)] If $b$ divides $a$, then there exists a
maximal solvable primitive subgroup of $\GL(e,p^b)$ with
generators $x_1,\ldots,x_l$, and matrices  $t\in \GL(a/b,p^b)$ of
order $p^a - 1$
and $s\in \GL(a/b,p^b)$ of order $a/b$ satisfying $t^s=t^{p^b}$,
such that the subgroup
$$\langle t\otimes I_e, s\otimes I_e,
I_{a/b} \otimes x_1,\ldots, I_{a/b} \otimes x_l\rangle
\leq \GL(d/b,p^b):\langle
s\otimes I_e\rangle\leq \GL(d,p)$$ is conjugate to a normal
subgroup of $H$ containing~$A$ and of index dividing~$b$.\smallskip

\item[(ii)] If $b$ does not divide $a$, then there exists a
maximal solvable primitive subgroup of $\GL(e,p)$ with
generators $x_1,\ldots,x_l$,
and matrices
$t\in \GL(a,p)$ of order $p^a-1$ and
$s\in \GL(a,p)$ of order $a$ satisfying $t^s=t^p$,
such that $H$ is conjugate to the subgroup
$$\langle t\otimes I_e,s\otimes I_e, I_a\otimes x_1,\ldots,
 I_a \otimes x_l
\rangle\leq\GL(d,p).$$
\enmrt
\elmm

\proof (i) By Lemma~\ref{211220b}(i), $F=U\circ E$, where $E$ is an extraspecial group of order $r^{2k+1}$; moreover, $F$ is unique up to conjugation in $\GL(d,p)$.
%
%
By \cite[2.5.14]{Short1992}, $U=\langle z\otimes I_e \rangle$ where
$z$ is a Singer cycle of $\GL(a,p)$.
By \cite[Theorem~2.5.15]{Short1992}, $L=C_{\GL(d,p)}(U)\cong \GL(e,p^a)$,
and $N_{\GL(d,p)}(U)=\GL(e,p^a):\langle \psi\rangle$,
where $\psi$ is a field automorphism of order $a$. Identifying $L$ with $\GL(e,p^a)$, we have the series of subgroups
$$1<U<F<A\le N_L(F)\le L=\GL(e,p^a).$$
Furthermore, $N_L(F)\cong\Sp(2k,r)$ (cf. Lemma~\ref{211220b}(ii)).\medskip

Since $b$ divides $a$, there is an embedding $\GL(a/b,p^b)\leq \GL(a,p)$ and we may choose it so that $z$ lies in its image. Let $t$ be the
preimage of $z$ under this embedding. By \cite[Lemma~2.7]{PTV2004},
there exists $s\in\GL(a/b,p^b)$ of order $a/b$ such that $t^s=t^{p^b}$.\medskip

The embeddings $\GF(p)\le \GF(p^b)\le \GF(p^a)$ yield embeddings
$$\GL(e,p^a)\le \GL(e\cdot (a/b),p^b)\le \GL(d,p).$$
In particular, $A=C_{H}(U)$ is a subgroup of $\GL(e\cdot (a/b),p^b)$.\medskip
%

By Lemma~\ref{211220b}(i), $L_1=\GL(e,p^b)$ contains a subgroup 
$F_1=U_1\circ E$, 
where $U_1=Z(L_1)$ and $E$ is extraspecial of order $r^{2k+1}$.
Hence we have the series of subgroups
$$1<U_1<F_1<A_1\le N_{L_1}(F_1)\le L_1=\GL(e,p^b).$$
It follows from Lemma~\ref{211220b}(ii) that $N_{L_1}(F_1)\cong\Sp(2k,r)\cong N_L(F)$. Since $U= \langle t\otimes I_e \rangle$ commutes with all matrices from $I_{a/b}\otimes N_{L_1}(F_1)$,
$$N_{\GL(e\cdot (a/b), p^b)}(F)=\langle U, I_{a/b}\otimes N_{L_1}(F_1)\rangle.$$ Therefore $A\le\langle U, I_{a/b}\otimes N_{L_1}(F_1)\rangle$ and
there exists $A_1=\langle x_1,\ldots,x_l\rangle\le L_1$ such that $A=U\circ (I_{a/b}\otimes A_1)$.

Since $H/A$ is cyclic of order dividing $a$ and $s\otimes I_e$ commutes with all matrices from $I_{a/b}\otimes A_1$,
it follows that $\langle A, s\otimes I_e\rangle$ is conjugate to a normal subgroup of $H$ of index dividing $b$, as required.\medskip

It remains to show that $A_1$ is primitive in $\GL(e,p^b)$. By way of contradiction, assume that there is a proper $\GF(p^b)$-subspace $\ov W$ of the natural $\GF(p^b)A_1$-module $\ov V$ with $\ov V=\bigoplus_{\ov g\in A_1}\ov W^g$. Consider an embedding of $\GL(e,p^b)$ into $\GL(e,p^a)$ such that $A=A_1\circ U$.
Therefore
$$
\wt V=\GF(p^a)\otimes_{\GF(p^b)}\ov{V}=\bigoplus_{g\in A}\wt{W}^g,
$$
where $\wt V$ is the natural $\GF(p^a)A$-module
and $\wt W=\GF(p^a)\otimes_{\GF(p^b)}\overline{W}$. Hence $A$ is imprimitive in $\GL(e,p^a)$.

Since $U=\langle z\otimes I_e \rangle$ where $z$ is a Singer cycle of $\GL(a,p)$, the $\GF(p^a)$-subspaces $\wt W^g$ are $U$-invariant for all $g\in A$. Therefore
$$
V=\GF(p^a)\otimes_{\GF(p)}\wt{V}=\bigoplus_{g\in A}{W}^g,
$$
where $V$ is the natural $\GF(p)A$-module and $W=\GF(p^a)\otimes_{\GF(p)}\wt{W}$.

Now $H=\langle A,x\rangle,$ where $x$ induces a field automorphism 
$\psi$ of $U$ (cf.\ Theorem~\ref{060519a}(iv)). 
Since $W^g$ is $\psi$-invariant for each $g\in A$,
$$V=\bigoplus_{h\in H}W^h,$$ which contradicts the primitivity of
$H$ in $\GL(d,p)$.\medskip

(ii) We put $b=1$ and
choose an extraspecial subgroup $F_1\le L_1=\GL(e,p)$
such that $N_{L_1}(F_1)\cong N_L(F)$
(cf.\ Lemma~\ref{211220b});
the choice depends on which of two conjugacy classes contains $F$.
The rest
of the proof is similar to that of~(i). \eprf

\section{Primitive solvable linear groups: computations}\label{090120h}

To complete the proof of Theorem~\ref{241219a}, it suffices to establish
the following.

\thrml{301219b}
Let $G$ be a linearly primitive solvable permutation group with
exceptional parameters listed in columns~$2$--$5$
of Table~$\ref{t2}$.  Then $G^{(3)}$ is solvable.
\ethrm

\proof Let $G_0$ be a point stabilizer of~$G$ with underlying vector
space $V$. By the monotonicity of
the $3$-closure operator, we may assume that $H=G_0$ is a maximal
solvable primitive linear group. 

Our proof is computational; more details are given in Section \ref{add-details}.
For each choice of the parameters $p$, $d$, $a$, $e$,
we compute a list $\cH=\cH(p,d,a,e)$  of
solvable primitive subgroups $H$ of $\GL(d, p)$,
ensuring that $\cH$ includes representatives of all conjugacy
classes of such {\it maximal} solvable subgroups
for specified $a$ and $e$.


For every $H\in\cH$, we search for nonzero $\alpha \in V$ such
that one of the following conditions is satisfied:
\nmrt
\tm{A} $\alpha^H$ is a regular orbit of $H$,
so Corollary~\ref{060519b} can be applied;
\tm{B} the restriction of $H$ to $\alpha^H$ is $2$-closed,
so Lemma \ref{220520c} can be applied.
\enmrt
Of course, (A) implies (B).
If we find such $\alpha$, then $G$ is $3$-closed
by Corollary~\ref{060519b} and Lemma~\ref{220520c}, respectively; so $G^{(3)}=G$ is solvable. For those
groups~$H$ where no such $\alpha$ is found,
we verify that $H$ acts transitively on the nonzero vectors of~$V$.
In this case, $G^{(3)}$ is solvable by Lemma~\ref{220520b}.


Our results are summarized in
Table~\ref{t2}.
The $7$th column lists the numbers of groups $H\in\cH$
for which conditions (A) or (B) are satisfied;
we indicate (A) or~(B) by writing ``partly regular'' or
``$2$-closed constituent'', respectively; we  write ``transitive'' if $H\in\cH$  acts transitively on the nonzero vectors of~$V$.  Detailed results, including the {\sf GAP} procedures used, generators of
$H$, vectors~$\alpha$, and certificates for (A) and (B),
are available at~\cite{Comp}.\eprf

\subsection{Constructing $\cH$}\label{add-details}
The construction of the list $\cH$ for given parameters $p$, $d$, $a$, $e$
naturally divides into two cases I and II; the $6$th column of Table \ref{t2} indicates which of them is applied to the stated parameters $p$, $d$, $a$, $e$.

\subsubsection*{\bf Case I}
Where possible, using the {\sf GAP} package IRREDSOL,
we constructed the list ${\cL}$ of all solvable primitive
subgroups of $\GL(d,p)$.
By Theorem~\ref{060519a}, every group $H\in\cH$ has
order $(p^a-1)\,e^2\,s(H)\,a'$, where $s(H)=|M|$ is the order of
a maximal solvable subgroup $M\cong A/F$ of the corresponding linear group
$S$ from Lemma~\ref{spo}, and $a'$ is a divisor of~$a$.  Table~\ref{t1} lists the
possible values of~$s(H)$ for all relevant values of~$e$.
By filtering ${\cL}$ with respect to the possible orders,
we obtain $\cH$.

\subsubsection*{\bf Case II}
We use auxiliary results from Section~\ref{280520c} and \cite[Section~2.5]{Short1992} together with computations
in {\sc Magma} to construct~$\cH$.
Recall that $e=r^k$ where $p\neq r\in\{2,3\}$, $d=ae$, and $b$ is the least positive integer with $p^b\equiv1\pmod{r^c}$ where $c=2$ for $r=2$ and $c=1$ otherwise.
Since $r\in\{2,3\}$, we deduce that $b\leq 2$.
For each remaining set of parameters, we proceed as follows.

\begin{enumerate}
\item[1.]
If $b$ divides $a$ and IRREDSOL contains the solvable primitive
subgroups of $\GL(e,p^b)$, or $b$ does not divide $a$ and IRREDSOL contains
the solvable primitive subgroups of $\GL(e,p)$, then
we construct a list $\cH_0$ by
extracting from the relevant output those groups containing
an extraspecial subgroup of order $r^{2k+1}$ and proceed to Step~5.

\item[2.] Otherwise, using Holt's implementation in
{\sc Magma} of the algorithm
of \cite{HoltCRD}, we construct in $\GL(e,p^b)$ extraspecial subgroups $E$ of order $r^{2k+1}$ (one if~$b$ divides $a$, and two if not) and the corresponding subgroups $F=E\circ U$, where $U=Z(\GL(e,p^b))$, and normalizers $N=N_{\GL(e,p^b)}(F)$, where $N/F\cong\Sp(2k,r)$ or $O^\varepsilon(2k,2)$ (cf. items (i) and (ii) of Lemma~\ref{211220b}).

\item[3.] Lemma~\ref{211220b}(iii) implies that $A/F$ is a maximal solvable subgroup of $N/F$ such that the $r$-radical of $A/F$
is trivial if $N/F\cong\Sp(2k,r)$ and otherwise has order at most $2$.
Using standard tools in {\sc Magma}, we construct the list $\cL$ consisting of all maximal solvable subgroups of $\Sp(2k,r)$ with trivial $r$-radical (if $b$ divides $a$), or of all maximal solvable subgroups of $O^+(2k,2)$
and $O^-(2k,2)$ with $r$-radical of order at most $2$ (if $b$ does not divide $a$).

\item[4.] Following \cite[Theorems~2.5.35 and~2.5.37]{Short1992}, for each subgroup in $\cL$, we produce generators of its complete preimage in $N$.
Thus we obtain the list $\cH_0$ containing up to conjugation all maximal solvable primitive
subgroups of $\GL(e,p^b)$ if $b$ divides $a$, and an equivalent list
in $\GL(e,p)$ if not.

\item[5.] Suppose $b$ divides $a$.
By Lemma \ref{Kron}(i), every maximal solvable primitive
subgroup $H$ of $\GL(d,p)$ contains up to conjugation 
an appropriate normal subgroup $H_1$ 
of index $b \leq 2$ in $H$.
We construct
$t,s$ as in Lemma \ref{Kron}(i).
For $H_0 \leq \GL(e,p^b)$ in $\cH_0$, define
$$H_1=\langle t\otimes I_e,s\otimes I_e,I_{a/b}\otimes H_0\rangle \leq\GL(d,p).$$
If $b=1$, then we take as $\cH$ the set consisting of $H := H_1$
for every $H_0 \in \cH_0$.
If $b=2$, then let $F=\Fit(H_1)$ and define
$$L=N_{\GL(d,p)}(F)\cong F:(\Sp(2k,r):\Z_a).$$
Now we take as $\cH$ the set
consisting of $H := N_L(H_1)$ for every $H_0 \in \cH_0$.

Suppose $b$ does not divide $a$.
We construct $t,s$ as in Lem\-ma~\ref{Kron}(ii),
and take as $\cH$ the set
$$\{\langle t\otimes I_e,s\otimes I_e,I_{a}\otimes
H_0\rangle\mid H_0\in\cH_0\}.$$ The same lemma guarantees that up to conjugation all maximal solvable primitive subgroups of $\GL(d,p)$ are in~$\cH$.

\end{enumerate}

\subsection{Processing $\cH$}
We discuss briefly how we process each $H \in \cH$.
For given $\alpha \in V$, condition (A) is readily checked. Since the {\tt TwoClosure} command in {\sf GAP} is time consuming,
we use the {\sf GAP} package COCO2P~\cite{KlinCOCO2P} to verify (B).
Namely, we compute the automorphism group~$\wt H$ of the coherent configuration associated with the restriction of $H$ to $\alpha^H$ 
using the COCO2P commands {\tt ColorGraph} and {\tt AutomorphismGroup}. Now $H$ is $2$-closed if and only if $H=\wt H$ \mbox{\cite[Corollary~2.2.18]{CP2019}}.

It was sometimes infeasible to compute all $H$-orbits on $V$,
so we randomly selected vectors $\alpha$
until we found one which satisfies either (A) or (B).
Hence, in principle, some cases resolved by (B) could also
be resolved by (A).

\medskip
{\small
\begin{longtable}{|r|r|r|r|r|r|l|l|}
\caption{Results for  groups with exceptional parameters\label{t2}} \\
\hline
No.  &  $e$      &  $p$    & $d$   &  $a$      &  case & results \\
\hline
$1$  & $16$      & $3$    &   $16$       &  $1$       &  II                          &$780$, partly regular  \\
         &               &           &                    &               &                                & $21$, $2$-closed constituent  \\
\hline
$2$  &  $16$     & $5$    &   $16$       &  $1$       & II                           &$1085$, partly regular    \\
\hline
\hline
\hline
$3$  &  $9$       & $2$    & $18$         &  $2$       & I                            &  $31$, $2$-closed constituent \\ 
\hline
$4$  &  $9$       & $7$     &   $9$         &  $1$      &  II                          & $44$, partly regular  \\
\hline
$5$  &  $9$       & $13$   & $9$           &  $1$      & II                          &  $44$, partly regular   \\
\hline
$6$  &  $9$       & $2$     &   $36$         &  $4$      &  II                          & 7, partly regular  \\
\hline
$7$  &  $9$     & $19$   & $9$                &  $1$      & II                         &   $44$, partly regular  \\
\hline
$8$  &  $9$       & $5$     &   $18$         &  $2$      &  II                          & $44$, partly regular  \\
\hline
\hline
\hline
$9$  &   $8$     & $3$    & $8$          &  $1$       & I                        &  $6$,  $2$-closed constituent\\
\hline
$10$  &   $8$     & $5$    & $8$           &  $1$     & I                         &  $4$, $2$-closed constituent   \\
\hline
$11$  &   $8$     & $7$    & $8$           &  $1$       & II                        & $30$, partly regular  \\
          &               &           &                    &               &                          &  $1$, $2$-closed constituent \\
\hline
$12$  &  $8$       & $3$     &   $16$         &  $2$      &  II                  & $63$, partly regular  \\
\hline
$13$  &   $8$    & $11$   & $8$           &  $1$       &  II                      & $122$, partly regular  \\
\hline
$14$  &   $8$    & $13$   & $8$           &  $1$       &  II                      &   $63$, partly regular\\
\hline
$15$  &   $8$    & $17$   & $8$           &  $1$       &  II                      &  $63$, partly regular \\
\hline
$16$  &   $8$    & $19$   & $8$           &  $1$       &  II                      &  $123$, partly regular \\
\hline
$17$  &   $8$    & $5$   & $16$           &  $2$       &  II                    &  $4$, partly regular\\
\hline
$18$  &   $8$    & $3$   & $24$           &  $3$       &  II                    &   $6$, partly regular \\
\hline
\hline
\hline
$19$  &   $4$    & $3$     & $4$          &  $1$       & I                         & $3$,  $2$-closed constituent \\
\hline
$20$   &   $4$    & $5$    & $4$           &  $1$      & I                        & $2$, $2$-closed constituent\\
\hline
$21$  &   $4$     & $7$    & $4$          &  $1$       & I                        & $13$,  $2$-closed constituent \\
\hline
$22$  &   $4$     & $3$    & $8$          &  $2$       & I                        &  $22$,  $2$-closed constituent \\
\hline
$23$  &   $4$     & $11$   & $4$         &  $1$       & I                        & $3$, $2$-closed constituent   \\
\hline
$24$  &   $4$     & $13$  & $4$          &  $1$      & I                        & $2$,  $2$-closed constituent   \\
\hline
$25$  &   $4$      & $17$  & $4$          &  $1$     & I                        & $2$, $2$-closed constituent   \\
\hline
$26$  &   $4$      & $19$  & $4$         &  $1$     & I                        & $7$, partly regular  \\
\hline
$27$  &   $4$      & $23$ & $4$          &  $1$     & I                        & $11$, partly regular \\
           &               &           &                  &              &                          &  $1$, $2$-closed constituent \\
\hline
$28$  &   $4$       & $5$  & $8$          &  $2$    & I                        &  $41$, partly regular \\
           &               &           &                  &              &                          &  $21$, $2$-closed constituent \\
 \hline
$29$  &   $4$       & $3$  & $12$       &  $3$    & I                        &  $9$, partly regular\\
\hline
$30$  &  $4$       & $29$ & $4$          &  $1$     & I                       & $2$, $2$-closed constituent  \\
\hline
$31$  &   $4$       & $31$  & $4$        &  $1$     & I                       & $16$, partly regular  \\
\hline
$32$  &   $4$      & $37$   & $4$       &  $1$     & I                        & $2$, $2$-closed constituent   \\
\hline
$33$  &   $4$       & $41$  & $4$        &  $1$    & I                        & $2$, $2$-closed constituent   \\
\hline
$34$  &   $4$       & $43$  & $4$       &  $1$    & I                        & $6$, $2$-closed constituent   \\
\hline
$35$  &   $4$       & $47$  & $4$       &  $1$    & I                        & $24$, partly regular  \\
           &               &           &                  &              &                          &  $2$, $2$-closed constituent  \\
\hline
$36$  &   $4$       & $7$     & $8$       &  $2$    & I                        &  $11$, partly regular   \\
           &               &           &                  &              &                          &  $17$, $2$-closed constituent\\
\hline
$37$  &   $4$       & $53$  & $4$       &  $1$    & I                       & $2$, $2$-closed constituent  \\
\hline
$38$ &   $4$         & $59$  & $4$       &  $1$    & I                      & $2$, $2$-closed constituent  \\
          &                   &            &               &            &                          &  $1$, partly regular \\
\hline
$39$ &   $4$         & $61$    & $4$       &  $1$   & I                     & $2$, $2$-closed constituent \\
\hline
$40$  &   $4$        & $67$   & $4$        &  $1$   &    II                    & $2$, partly regular \\
\hline
$41$  &   $4$        & $71$    & $4$       &  $1$     &  II                  & $2$, partly regular  \\
\hline
$42$  &   $4$        & $3$      & $16$       &  $4$    &  II                   & $2$, partly regular  \\
\hline
$43$  &   $4$        & $11$      & $8$       &  $2$    &  II                   &  $2$, partly regular \\
\hline
$44$  &   $4$        & $5$      & $12$       &  $3$    &  II                   &  $2$, partly regular  \\
\hline
$45$  &   $4$        & $13$      & $8$       &  $2$    &  II                   &  $2$, partly regular\\
\hline
$46$  &   $4$        & $3$      & $20$       &  $5$    &  II                   &  $3$, partly regular  \\
\hline
\hline
\hline
$47$  &   $3$         & $2$     & $6$        &  $2$     & I                  & $2$,  $2$-closed constituent   \\
\hline
$48$  &   $3$         & $7$     & $3$        &  $1$     & I                  & $1$,  $2$-closed constituent    \\
\hline
$49$  &   $3$         & $13$   & $3$        &  $1$     & I                  & $1$,  $2$-closed constituent   \\
\hline
$50$  &   $3$         & $2$   & $12$        &  $4$     & I                  &  $3$,  $2$-closed constituent   \\
\hline
$51$  &   $3$        & $19$    & $3$        &  $1$     & I                  & $1$,  $2$-closed constituent \\
\hline
$52$ &   $3$         & $5$      & $6$        &  $2$     & I                  & $4$,  $2$-closed constituent \\
\hline
$53$ &   $3$         & $7$       & $6$        &  $2$     & I                  & $5$,  partly regular  \\
\hline
$54$ &   $3$         & $2$       & $18$      &  $6$     & I                  & $7$,  partly regular  \\
&&&&&I&$11$, $2$-closed constituent    \\
\hline
$55$ &   $3$         & $11$      & $6$        &  $2$     & I                  & $4$,  partly regular   \\
\hline
$56$ &   $3$         & $13$     & $6$        &  $2$     & I                  &  $4$,  partly regular    \\
\hline
$57$ &   $3$       & $2$     & $24$      &  $8$     & II                & $1$, partly regular \\
\hline
$58$ &   $3$       & $17$    & $6$        &  $2$     & II                & $1$, partly regular   \\
\hline
$59$ &   $3$       & $7$      & $9$         &  $3$     & II               & $1$, partly regular  \\
\hline
$60$ &   $3$       & $19$    & $6$         &  $2$     & II               & $1$,  partly regular  \\
\hline
\hline
\hline
$61$   &   $2$      & $3$      & $2$        &  $1$      & I                 & $2$, transitive  \\
\hline
$62$   &   $2$      & $5$      & $2$        &  $1$     & I                  &  $1$, transitive \\
\hline
$63$  &   $2$       & $7$      & $2$       &  $1$      & I                   & $5$, $2$-closed constituent  \\
\hline
$64$  &   $2$       & $3$     & $4$       &  $2$      & I                   & $6$, $2$-closed constituent\\
\hline
$65$  &   $2$       & $11$    & $2$       &  $1$       & I                  & $3$, $2$-closed constituent \\
\hline
$66$  &   $2$       & $13$   & $2$        &  $1$       & I                  & $2$, $2$-closed constituent \\
\hline
$67$  &   $2$       & $17$   & $2$         &  $1$       & I                  & $1$, $2$-closed constituent \\
\hline
$68$  &   $2$       & $19$    & $2$        &  $1$       & I                & $3$,  $2$-closed constituent  \\
\hline
$69$  &   $2$       & $23$   & $2$        &  $1$        & I                & $7$,  $2$-closed constituent  \\
\hline
$70$  &   $2$       & $5$      & $4$        &  $2$        & I               & $20$,  $2$-closed constituent  \\
\hline
$71$  &   $2$       & $3$      & $6$        &  $3$        & I               & $3$,  partly regular   \\
           &                 &              &              &                 &                  & $1$, $2$-closed constituent   \\
\hline
$72$  &   $2$       & $29$    & $2$         &  $1$       & I               & $1$,  $2$-closed constituent    \\
\hline
$73$  &   $2$       & $7$       & $4$         &  $2$       & I              & $12$,  $2$-closed constituent    \\
\hline
$74$  &   $2$       & $3$       & $8$         &  $4$       & I              & $30$,  $2$-closed constituent   \\
\hline
$75$  &   $2$       & $11$     & $4$         &  $2$       & I               & $14$,  $2$-closed constituent      \\
\hline
$76$  &   $2$       & $5$       & $6$         &  $3$       & I               & $2$, $2$-closed constituent  \\
\hline
$77$  &   $2$       & $13$     & $4$         &  $2$       & I               & $8$,  partly regular   \\
\hline
$78$  &   $2$       & $3$       & $10$       &  $5$       & I               & $4$,  partly regular   \\
\hline
$79$  &   $2$       & $17$      & $4$         &  $2$       & I               & $8$,  $2$-closed constituent   \\
          \hline
$80$  &   $2$       & $7$        & $6$         &  $3$       & I               & $22$,  partly regular   \\
\hline
$81$  &   $2$       & $19$      & $4$         &  $2$       & I               & $19$,  $2$-closed constituent   \\
           \hline
$82$  &   $2$       & $23$      & $4$         &  $2$       & I               & $12$,  partly regular \\
\hline
$83$  &   $2$       & $5$        & $8$         &  $4$       & I               & $27$,  partly regular   \\
\hline
$84$  &   $2$       & $3$        & $12$       &  $6$       & I                & $17$,  partly regular     \\
\hline
$85$  &   $2$       & $29$      & $4$         &  $2$       & I                 & $8$,  partly regular     \\
\hline
$86$  &   $2$       & $31$      & $4$         &  $2$       & I                 & $8$,  partly regular    \\
\hline
$87$  &   $2$       & $11$      & $6$         &  $3$       & I                 & $7$,  partly regular    \\
\hline
$88$  &   $2$       & $37$      & $4$         &  $2$       & I                & $8$,  partly regular    \\
\hline
$89$  &   $2$       & $41$      & $4$         &  $2$       & I                & $8$,  partly regular    \\
\hline
$90$  &   $2$       & $43$      & $4$         &  $2$       & I                & $8$,  partly regular    \\
\hline
$91$  &   $2$       & $3$        & $14$        &  $7$       & I                & $4$,  partly regular    \\
\hline
$92$  &   $2$       & $13$       & $6$         &  $3$       & I                & $5$,  partly regular    \\
\hline
$93$  &   $2$       & $47$       & $4$        &  $2$       & I                  & $8$,  partly regular   \\
      &             &            &            &            & I               & $1$,  $2$-closed constituent   \\
\hline
$94$  &   $2$       & $7$         & $8$         &  $4$       & I                 & $23$,  partly regular      \\
 \hline
$95$  &   $2$       & $53$       & $4$         &  $2$       & I                 & $8$,  partly regular    \\
\hline
$96$  &   $2$       & $5$         & $10$       &  $5$       & I                 & $2$,  partly regular   \\
\hline
$97$  &   $2$       & $59$       & $4$         &  $2$       & I                  & $8$,  partly regular   \\
\hline
$98$  &   $2$       & $61$       & $4$         &  $2$       & I                  & $8$,  partly regular   \\
\hline
$99$  &   $2$       & $67$       & $4$         &  $2$       &   II                &  $1$, partly regular \\
\hline
$100$  &   $2$       & $17$       & $6$         &  $3$      &   II              & $1$, partly regular   \\
\hline
$101$  &   $2$       & $71$       & $4$         &  $2$       &  II                & $4$, partly regular  \\
\hline
$102$  &   $2$       & $73$       & $4$         &  $2$       & II                & $1$, partly regular     \\
\hline
\end{longtable}
}

\end{document}